\documentclass[11pt]{article}

\usepackage{amssymb}
\usepackage{amsmath}
\usepackage{epsfig}
\usepackage{color}

\newtheorem{theorem}{Theorem}[section]
\newtheorem{lemma}[theorem]{Lemma}

\newtheorem{example}[theorem]{Example}

\renewcommand{\baselinestretch}{1.5}

\title{Almost $2$-perfect $8$-cycle systems}

\author{{\bf{Selda K\"{u}\c{c}\"{u}k\c{c}if\c{c}i} }\\ Department of
  Mathematics, Ko\c{c} University, Istanbul,  Turkey\\
   email: {\tt skucukcifci@ku.edu.tr}\\
  {\bf{Charles Curtis Lindner}}\\
Department of Mathematics and Statistics, Auburn University, USA\\
email: {\tt lindncc@auburn.edu}\\
  {\bf{Sibel \"{O}zkan}}\\
Department of Mathematics, Gebze Technical University, Gebze, Kocaeli, Turkey\\
email: {\tt s.ozkan@gtu.edu.tr}\\
{\bf{Emine \c{S}ule Yaz\i c\i} }\\ Department of
  Mathematics, Ko\c{c} University, Istanbul,  Turkey\\
    email: {\tt eyazici@ku.edu.tr}}

\date{ }

\begin{document}
\maketitle\thispagestyle{empty}

\def\baselinestretch{1.5}\small\normalsize
\begin{abstract}
For an $m$-cycle $C$, an inside $m$-cycle of $C$ is a cycle on the same vertex set, that is edge-disjoint from $C$. In an $m$-cycle system,  $(\mathcal{X}, \mathcal{C})$, if inside $m$-cycles can be chosen -one for each cycle- to form another $m$-cycle system, then $(\mathcal{X}, \mathcal{C})$ is called an almost $2$-perfect $m$-cycle system. Almost $2$-perfect cycle systems can be considered as generalisations of $2$-perfect cycle systems. Cycle packings are generalisations of cycle systems that allow to have leaves after decomposition. In this paper, we prove that an almost $2$-perfect maximum packing of $K_n$ with $8$-cycles of order $n$ exists for each $n\geq 8$. We also construct a maximum $8$-cycle packing of order $n$ which is not almost $2$-perfect for each $n \geq 10$.
\end{abstract}

\noindent {\small AMS classification: $05B30$, $05B40$, $05C38$.}

\section{Introduction}

One of the oldest graph decomposition problems involves decomposing the complete graph $K_n$ into edge disjoint cycles. If all the cycles have uniform length, say $m$, then this decomposition is called an {\em $m$-cycle system}. More formally, we denote an $m$-cycle system of order $n$ by a pair $(\mathcal{X}, \mathcal{C})$ where $\mathcal{X}$ is an $n$-element set and $\mathcal{C}$ is a collection of edge-disjoint $m$-cycles which partitions the edge-set of $K_n$ with vertex set $\mathcal{X}$. To have an $m$-cycle system of order $n$, we need $n$ to be odd and $m$ to divide the total number of edges. Namely;


(i) $n \geq m \geq 3$,

(ii) $n\equiv 1$ (mod $2$), and

(iii) $n(n-1)/2m$ is an integer.

These conditions are shown to be sufficient when $n$ and $m$ have the same parity in \cite{AG}, and when $n$ and $m$ have different parity  in \cite{S}. Given $m$, the set of all $n$ satisfying these conditions is called the {\em spectrum} of the $m$-cycle system.

For a cycle $C$, let $c(i)$ be the set of edges that connects the vertices that are distance $i$ apart in $C$. For example, if $C$ is a 6-cycle, then $c(2)$ consists of two 3-cycles. If $C$ is a 5-cycle, then $c(2)$ is another 5-cycle. Given an $m$-cycle system $(\mathcal{X}, \mathcal{C})$, if the collection of $c(i)$ related to each $C$ forms another cycle system (not necessarily an $m$-cycle system), then $(\mathcal{X}, \mathcal{C})$ is called an {\em i-perfect m-cycle system.} There are many results on 2-perfect $m$-cycle systems where every pair of vertices is connected by a path of length two in exactly one $m$-cycle in $\mathcal{C}$. The interested reader is referred to \cite{AB, BL, L, LPR, LR} for results on $2$-perfect cycle systems.

Another interesting problem involving $m$-cycle systems is to replace each $m$-cycle $C$ in the system with another $m$-cycle $C'$ related with $C$, and then to see if the collection of the new $m$-cycles still gives an $m$-cycle system or not. We can for example relate $C'$ with $C$ by choosing it on the same vertex set. Given an $m$-cycle $C$, an edge-disjoint $m$-cycle $C'$ on the same vertex set is called an {\em inside $m$-cycle} of $C$.


If $C=(x_1,x_2,x_3,x_4,x_5)$ is a $5$-cycle, then there is only one inside $5$-cycle of $C$, namely $(x_1,x_3,x_5,x_2,x_4)$. So for $m=5$; being $2$-perfect is equivalent to inside $m$-cycles forming an $m$-cycle system. When $m$ is even, the cycles formed by $c(2)$ in the 2-perfect $m$-cycle system can not be $m$-cycles, as given in the example of a 6-cycle. That's why the generalisation of the idea as ``almost" 2-perfect cycle systems has gained attention lately.
For $m$-cycles when $m>5$ there are increasing number of possible inside cycles for a given cycle. A $6$-cycle $(x_1,x_2,x_3,x_4,x_5,x_6)$ has three inside $6$-cycles $(x_1,x_3,x_5,x_2,x_6,x_4)$, $(x_1,x_3,x_6,x_4,x_2,x_5)$ and $(x_1,x_4,x_2,x_6,x_3,x_5)$ for example. In \cite{LMR}, a $6$-cycle system $(\mathcal{X}, \mathcal{C})$ is called an {\em almost $2$-perfect $6$-cycle system} if it is possible to choose an inside $6$-cycle from each $6$-cycle in $\mathcal{C}$, so that the resulting collection of $6$-cycles is also a $6$-cycle system. They found the spectrum of almost $2$-perfect $6$-cycle systems.

A {\em cycle packing} of the complete graph $K_n$ is a triple $(\mathcal{X}, \mathcal{C},\mathcal{L})$ where $\mathcal{X}$ is the vertex set, $\mathcal{C}$ is a collection of edge-disjoint cycles from $K_n$, and the leave $\mathcal{L}$ is the collection of the edges in $K_n$ not belonging to any of the cycles in  $\mathcal{C}$. When $|\mathcal{L}|$ is smallest possible, then $(\mathcal{X}, \mathcal{C},\mathcal{L})$ is called a maximum packing. If $n$ is not in the spectrum of the $m$-cycle system, i.e; when it is not possible to decompose $K_n$ completely into $m$-cycles, it is wise to work on maximum packings of $m$-cycles of $K_n$.

In \cite{LMR}, Lindner et al. also proved the existence of almost $2$-perfect maximum packing of $K_n$ with $6$-cycles for each admissible $n$. In \cite{LM}, Lindner and Meszka considered the existence of almost $2$-perfect minimum covering of $K_n$ with $6$-cycles.

In this paper we will consider the next case; almost $2$-perfect $8$-cycle systems. Although a given $6$-cycle has the three possible inside $6$-cycles, an $8$-cycle has $177$ possible inside $8$-cycles. Let $(\mathcal{X}, \mathcal{C})$ be an $8$-cycle system, and let $\mathcal{C}'$ be a collection of inside $8$-cycles one from each of the cycles in $\mathcal{C}$. If $(\mathcal{X}, \mathcal{C}')$ is an $8$-cycle system, we will say that $(\mathcal{X}, \mathcal{C})$ is almost $2$-perfect. We show that the spectrum of almost $2$-perfect 8-cycle systems is same as the spectrum of 8-cycle systems, namely $n \equiv 1 \pmod{16}$.

For all other $n$ not in the spectrum of 8-cycle systems, we show the existence of almost $2$-perfect maximum packings with $8$-cycles.


Another question is: Given an $8$-cycle maximum packing $(\mathcal{X}, \mathcal{C},\mathcal{L})$, is it always possible to choose an inside $8$-cycle for each $8$-cycle in $\mathcal{C}$ so that the resulting collection of inside $8$-cycles is an $8$-cycle maximum packing? In other words, are all the 8-cycle maximum packings almost 2-perfect? The answer to this question is no, except for the orders $n=8$ and $n=9$. In the third section, we construct almost $2$-perfect maximum packings of $K_n$ with $8$-cycles for all $n\geq 8$, therefore we show that they exists for all admissible $n$. And in the fourth section, we construct maximum packings of $K_n$ with $8$-cycles that are not almost $2$-perfect for all $n\geq 10$ and present our observations from a comprehensive computer search.

\section{Preliminary results}

We will start by introducing the results that will be used throughout the paper. From now on, almost $2$-perfect will be abbreviated as A2P for brevity.

\begin{lemma}\label{K4t4s}
There exists an $A2P$ $8$-cycle decomposition of $K_{4t,4s}$, for all $t,s\in \mathbb{Z}^+$.
\end{lemma}

\noindent {\bf Proof} Let $X=\{x_0,x_1,x_2,x_3\}$ and $Y=\{y_0,y_1,y_2,y_3\}$ be the parts of $K_{4,4}$. Consider\\
$\mathcal{C}=\{(x_0,y_0, x_1,y_1,x_2,y_2,x_3,y_3)$, $(x_0,y_2,x_1,y_3, x_2,y_0, x_3, y_1)\}$, and the inside cycles as\\
$\mathcal{C'}=\{(x_1,y_2,x_0,y_1,x_3,y_0,x_2,y_3)$, $(x_1,y_1,x_2,y_2,x_3,y_3,x_0,y_0)\}$ to get an A2P $8$-cycle decomposition of $K_{4,4}$.

Next let $\mathcal{X}=\{x_0,x_1,...,x_{4t-1}\}$ and $\mathcal{Y}=\{y_0,y_1,...,y_{4s-1}\}$ be the parts of $K_{4t,4s}$, where $X_i=\{x_{4i},x_{4i+1},x_{4i+2},x_{4i+3}\}$, $Y_j=\{y_{4j},y_{4j+1},y_{4j+2},y_{4j+3}\}$ for $i=0,1,2,...,t-1$, $j=0,1,...,s-1$ and $t,s\in \mathbb{Z}^+$. Placing an A2P $8$-cycle decomposition of $K_{4,4}$ on the vertex set $X_i \cup Y_j$ for each pair $i$, $j$ gives us an A2P $8$-cycle decomposition of $K_{4t,4s}$. \hfill $\square$

\begin{lemma} \label{ts}
There exists an $A2P$ $8$-cycle decomposition of $K_{4t,4s+2}$, for all $t, s \in \mathbb{Z}^+$.
\end{lemma}

\noindent {\bf Proof} Let $X=\{x_0, x_1,x_2,x_3\}$ and $Y=\{y_0, y_1,y_2,y_3,y_4,y_5\}$ be the parts of $K_{4,6}$. Consider \\
$\mathcal{C}=\{(x_0,y_4, x_1,y_1,x_2,y_2,x_3,y_3)$, $(x_0, y_0, x_1,y_3,x_2,y_4,x_3,y_5)$, $(x_0,y_1,x_3,y_0,x_2,y_5, x_1,y_2)\}$ and\\
$\mathcal{C'}=\{(x_1,y_2,x_0,y_1,x_3,y_4,x_2,y_3)$, $(x_1,y_4,x_0,y_3,x_3,y_0,x_2,y_5)$, $(x_1,y_1,x_2,y_2,x_3,y_5,x_0,y_0)\}$ to get an A2P $8$-cycle decomposition of $K_{4,6}$.

Next, let $\mathcal{X}=\{x_0,x_1,...,x_{4t-1}\}$ and $\mathcal{Y}=\{y_0,y_1,...,y_{4s+1}\}$ be the parts of $K_{4t,4s+2}$, where $X_i=\{x_{4i},x_{4i+1},x_{4i+2},x_{4i+3}\}$, $Y_0=\{y_0,y_1,y_2,y_3,y_4,y_5\}$, $Y_j=\{y_{4j+2},y_{4j+3},y_{4j+4},y_{4j+5}\}$ for $i=0,1,...,t-1$, $j=1,2,...,s-1$ and $t,s\in \mathbb{Z}^+$. Placing an A2P $8$-cycle decomposition of $K_{4,4}$ on the vertex set $X_i \cup Y_j$ for each pair $i$, $j$, where $i\geq 0$, $j\geq 1$ and an A2P $8$-cycle decomposition of $K_{4,6}$ on the vertex set $X_i \cup Y_0$ for each $i=0,1,...,t-1$ gives us an A2P $8$-cycle decomposition of $K_{4t,4s+2}$. \hfill $\square$

\begin{lemma}\label{hole}
If there exist an $A2P$ $8$-cycle system of order $r+1$ and an $A2P$ $8$-cycle decomposition of $K_{r,s}$, then there exists an $A2P$ $8$-cycle decomposition of $K_{r+s+1}\setminus K_{s+1}$.
\end{lemma}

\noindent {\bf Proof} Let $\mathcal{X}=\{\infty\}\cup \{x_1,x_2,...,x_{r}\}\cup \{y_1,y_2,...,y_{s}\}$. Placing an A2P $8$-cycle system of order $r+1$ on $\{\infty\}\cup \{x_1,x_2,...,x_{r}\}$ and an A2P $8$-cycle decomposition of $K_{r,s}$ on $\{x_1,x_2,...,x_{r}\}\cup \{y_1,y_2,...,y_{s}\}$ gives an A2P maximum $8$-cycle decomposition of $K_{r+s+1}\setminus K_{s+1}$, where the vertex set of $K_{s+1}$ is $\{\infty\}\cup \{y_1,y_2,...,y_{s}\}$.  \hfill $\square$

\begin{lemma}\label{bipar}
If there exist an $A2P$ maximum $8$-cycle packing of order $r$ and of order $s$ with a $1$-factor leave and an $A2P$ $8$-cycle decomposition of $K_{r,s}$, then there exists an $A2P$ maximum $8$-cycle packing of order $r+s$ with a $1$-factor leave.
\end{lemma}

\noindent {\bf Proof} Let $\mathcal{X}=\{x_1,x_2,...,x_{r}\}\cup \{y_1,y_2,...,y_{s}\}$. Placing an A2P maximum $8$-cycle packing of order $r$ on $\{x_1,x_2,...,x_{r}\}$, an A2P maximum $8$-cycle packing of order $s$ on $\{y_1,y_2,...,y_{s}\}$, and an A2P $8$-cycle decomposition of $K_{r,s}$ on $\{x_1,x_2,...,x_{r}\}\cup \{y_1,y_2,...,y_{s}\}$ gives an A2P maximum $8$-cycle packing of order $r+s$.  \hfill $\square$\vspace{0.2in}

\begin{lemma}\label{biparhole}
If there exist an $A2P$ maximum $8$-cycle packing of order $r$ and an $A2P$ $8$-cycle decomposition of $K_{r,s}$, then there exists an $A2P$ maximum $8$-cycle packing of $K_{r+s}\setminus K_{s}$.
\end{lemma}

\noindent {\bf Proof} Let $\mathcal{X}=\{x_1,x_2,...,x_{r}\}\cup \{y_1,y_2,...,y_{s}\}$. Placing an A2P maximum $8$-cycle packing of order $r$ on $\{x_1,x_2,...,x_{r}\}$ and an A2P $8$-cycle decomposition of $K_{r,s}$ on $\{x_1,x_2,...,x_{r}\}\cup \{y_1,y_2,...,y_{s}\}$ gives an A2P maximum $8$-cycle packing of $K_{r}\setminus K_{s}$.  \hfill $\square$\vspace{0.2in}

Now we can present the main construction that will be used to construct almost 2-perfect maximum 8-cycle packings.

\noindent{\bf Main Construction}

Let $H$ be a finite set and $k$ be a positive integer.
Then, let $X=H\cup \{(i,j)\ | \ 1\leq i\leq k, 1\leq j\leq 16\}$.

\noindent (1) On $H\cup \{(1,j)\ | \ 1\leq j\leq 16\}$, place an A2P maximum $8$-cycle packing of order $16+h$, where $h=|H|$.

\noindent (2) On each set $H\cup \{(i,j)\ | \
1\leq j\leq 16\}$, for $2\leq i\leq k$, place an A2P $8$-cycle decomposition of $K_{16+h}\setminus K_{h}$.

\noindent (3) For each $x, y= 1,2,...,k$, $x<y$ place an A2P $8$-cycle decomposition of $K_{16,16}$ on $\{(x,j)\ | \ 1\leq j\leq 16\} \cup \{(y,j)\ | \ 1\leq j\leq 16\}$.

Combining (1), (2), and (3) gives an A2P maximum $8$-cycle packing of order $16k+h$. \hfill $\square$

\section{Almost $2$-perfect maximum packings with $8$-cycles}

The following table gives leaves we considered for the maximum packings with $8$-cycles (see \cite{KLH} and \cite{LHG}).

\begin{center}
\begin{tabular}{|c|c|}
\hline
Spectrum for maximum packing with $8$-cycles  & Leave\\\hline
1 (mod 16)  & $\varnothing$\\ \hline
3 (mod 16)  & $C_3$\\ \hline
5 (mod 16)  & $K_5$\\ \hline
7 (mod 16)  & $C_5$\\ \hline
9 (mod 16), $n \neq 9$ & $C_4$\\ \hline
11 (mod 16)  & $C_3 + C_4$\\ \hline
13 (mod 16)  & bowtie\\ \hline
15 (mod 16)  & $C_4 + C_5
$\\ \hline
0, 2, 8, 10 (mod 16)  & $1$-factor\\ \hline
4, 6, 12, 14 (mod 16)  & $K_4 \ +$  a $1$-factor on the remaining vertices\\ \hline
\end{tabular}\vspace{0.2in}

Table 1: Maximum packings with $8$-cycles
\end{center}

Now, let's analyse the cases given in Table 1 to see if there exists an A2P maximum packing for each of those. In the first case there is no leave, hence we have an 8-cycle system. For the rest of the cases we will work on maximum packings. In each case, we will start by giving small examples, then we will use these small cases in the main construction given above.

\begin{center} \underline{\large $n\equiv 1$ (mod $16$)} \end{center}

\begin{example} \label{ex17}
An $A2P$ $8$-cycle system of order $17$ exists.
\end{example}

\noindent  To show that, let $\mathcal{X}=\mathbb{Z}_{17}$ and $\mathcal{C}$ be the set of $8$-cycles with base $8$-cycle $(0,16,2,6,11,1,3,9)$ developed modulo $17$. The inside cycle $(0,2,1,6,9,16,3,11)$ of the base cycle of $\mathcal{C}$ forms a base $8$-cycle for $\mathcal{C}'$.

\begin{lemma} \label{1mod16}
For every $n\equiv 1 \ ($mod $16)$ with $n\geq 17$, there exists an $A2P$ $8$-cycle system of order $n$.
\end{lemma}

\noindent {\bf Proof} Let $h=1$ in the Main Construction. Since there exist an A2P $8$-cycle system of order $17$ given by Example \ref{ex17} and an A2P $8$-cycle decomposition of $K_{16,16}$ by Lemma \ref{K4t4s} the result follows.\hfill $\square$


\begin{center} \underline{\large $n\equiv 3$ (mod $16$)} \end{center}

\begin{example} \label{ex8}
An $A2P$ maximum $8$-cycle packing of order $8$ exists.
\end{example}

\noindent To show that, let $\mathcal{X}=\mathbb{Z}_8$ and consider the maximum 8-cycle packing \\
$\mathcal{C}=\{(1,4,3,6,5,2,0,7)$, $(1,3,2,4,7,5,0,6)$, $(1,2,7,6,4,5,3,0)\}$ with the leave \\
$\mathcal{L}=\{\{1,5\},\{3,7\},\{2,6\},\{4,0\}\}$.
Then the inside cycles \\
$\mathcal{C'}=\{(1,3,2,4,7,5,0,6)$, $(1,0,3,5,4,6,7,2)$, $(1,7,0,2,5,6,3,4)\}$ with the same leave forms another maximum packing. Hence, $(\mathcal{X}, \mathcal{C}, \mathcal{L})$ is an A2P maximum $8$-cycle packing.

\begin{example} \label{ex11}
An $A2P$ maximum $8$-cycle packing of order $11$ exists.
\end{example}

\noindent
Similarly, let $\mathcal{X}=\mathbb{Z}_{11}$ and consider \\
$\mathcal{C}=\{(0,3,1,4,2,5,7,6)$, $(0,4,6,1,5,8,9,10)$, $(0,5,3,2,8,10,7,9)$, $(0,7,1,9,2,10,3,8)$, \\
$(1,8,6,9,3,7,4,10)$, $(2,6,10,5,9,4,8,7)\}$ with the leave $\mathcal{L}=\{(0,1,2),(3,4,5,6)\}$. \\
Then we can choose the following inside cycles: \\
$\mathcal{C'}=\{(0,4,6,2,3,5,1,7)$, $(0,6,8,4,10,1,9,5)$, $(0,3,7,5,10,2,9,8)$, $(0,9,3,1,8,2,7,10)$, \\
$(1,6,10,3,8,7,9,4)$, $(2,5,8,10,9,6,7,4)\}.$

\begin{example} \label{ex19}
An $A2P$ maximum $8$-cycle packing of order $19$ exists.
\end{example}

\noindent Let $\mathcal{X}=\{\infty_1, \infty_2,\infty_3\}\cup \{x_1,x_2,...,x_8\}\cup \{y_1,y_2,...,y_8\}$. Place a copy of an A2P maximum $8$-cycle packing of order $11$ on $\{\infty_1, \infty_2,\infty_3\}\cup \{x_1,x_2,...,x_8\}$ and on $\{\infty_1, \infty_2,\infty_3\}\cup \{y_1,y_2,...,y_8\}$, where the $3$-cycle in the leaves is $(\infty_1, \infty_2,\infty_3)$ and the $4$-cycles are $(x_1,x_2,x_3,$ $x_4)$ and $(y_1,y_2,y_3,y_4)$, respectively. Then place an A2P maximum $8$-cycle packing of order $8$ on $\{x_1,x_2,x_3,x_4\} \cup \{y_1,y_2,y_3,y_4\}$ with the $1$-factor leave $\{\{x_1,x_3\}, \{x_2,x_4\}$, $ \{y_1,y_3\}, \{y_2,y_4\}\}$, and place A2P $8$-cycle decompositions of $K_{4,4}$ on $\{x_1,x_2,x_3,x_4\} \cup \{y_5,y_6,y_7,y_8\}$, on $\{x_5,x_6,x_7,$ $x_8\} \cup \{y_1,y_2,y_3,y_4\}$ and on $\{x_5,x_6,x_7,$ $x_8\} \cup \{y_5,y_6,y_7,y_8\}$.

\begin{lemma} \label{3mod16}
For every $n\equiv 3 \ ($mod $16)$ with $n \geq 19$, there exists an $A2P$ maximum $8$-cycle packing of order $n$.
\end{lemma}

\noindent {\bf Proof} There exist an A2P maximum $8$-cycle packing of order $19$ with a $3$-cycle leave given by Example \ref{ex19} which is also an A2P $8$-cycle decomposition of $K_{19}\setminus K_{3}$. An A2P $8$-cycle decomposition of $K_{16,16}$ also exists by Lemma \ref{K4t4s} as before, and the result follows by the Main Construction considering $h=3$. \hfill $\square$

\begin{center} \underline{\large $n\equiv 5$ (mod $16$)} \end{center}

\begin{example} \label{ex21}
An $A2P$ maximum $8$-cycle packing of order $21$ exists.
\end{example}

\noindent Considering $r=16$ and $s=4$ in Lemma \ref{hole} gives an A2P $8$-cycle decomposition of $K_{21}\setminus K_{5}$, which is an A2P maximum $8$-cycle packing of order $21$ with a $K_5$ leave on $\{\infty\}\cup \{y_1,y_2,y_3,y_4\}$ we were looking for.

\begin{lemma} \label{5mod16}
For every $n\equiv 5 \ ($mod $16)$ with $n \geq 21$, there exists an $A2P$ maximum $8$-cycle packing of order $n$.
\end{lemma}

\noindent {\bf Proof} There exist an A2P maximum $8$-cycle packing of order $21$ with a $K_5$ leave given in Example \ref{ex21}. Then the result follows by the Main Construction considering $h=5$. \hfill $\square$

\begin{center} \underline{\large $n\equiv 7$ (mod $16$)} \end{center}

\begin{example} \label{ex9}
An $A2P$ maximum $8$-cycle packing of order $9$ exists.
\end{example}

\noindent To see this, let $\mathcal{X}=\mathbb{Z}_{9}$, and consider \\
$\mathcal{C}=\{(0,2,4,1,5,6,7,8)$, $(0,4,3,1,6,8,5,7)$, $(0,5,2,7,3,8,4,6)$, $(1,7,4,5,3,6,2,8)$,\\
with the leave $\mathcal{L}=\{(0,1,2,3)\}$. Then choose the inside cycles as \\
$\mathcal{C'}=\{(0,4,5,2,8,6,1,7)$, $(0,6,4,1,8,7,3,5)$, $(0,2,6,3,4,7,5,8)$, $(1,5,6,7,2,4,8,3)$.

\begin{example} \label{ex15}
An $A2P$ maximum $8$-cycle packing of order $15$ exists.
\end{example}

\noindent To see this, let $\mathcal{X}=\mathbb{Z}_{15}$ and consider the maximum packing  \\
$\mathcal{C}=\{(0,1,8,14,4,5,6,12)$, $(0,2,11,3,10,13,4,9)$, $(0,3,6,1,5,12,7,4)$, $(0,5,8,13,7,3,9,6)$,\\
$(0,8,2,6,13,9,1,10)$, $(0,11,14,12,2,9,5,13)$, $(9,12,3,1,4,10,2,14)$, $(1,13,2,4,8,10,5,14)$, \\
$(0,7,1,11,4,6,10,14)$, $(1,2,3,5,7,11,10,12)$, $(2,5,11,13,3,8,6,7)$, $(3,4,12,8,11,9,7,14)$,\\
with the leave $\mathcal{L}=\{(6,11,12,13,14),(7,8,9,10)\}$. \\
Now we can choose the inside cycles as\\
$\mathcal{C'}=\{(0,4,6,8,12,1,14,5)$, $(0,11,4,10,2,9,13,3)$, $(0,1,3,4,12,6,5,7)$, $(0,8,3,5,9,7,6,13)$,\\
$(0,2,1,13,8,10,6,9)$, $(0,14,9,11,13,2,5,12)$, $(9,1,10,12,14,3,2,4)$, $(1,8,2,14,10,13,4,5)$,\\
$(0,6,1,4,14,7,11,10)$, $(1,7,2,12,3,10,5,11)$, $(2,11,8,5,13,7,3,6)$, $(3,11,14,8,4,7,12,9)$.

\begin{example} \label{ex23}
An $A2P$ maximum $8$-cycle packing of order $23$ exists.
\end{example}

\noindent Let $\mathcal{X}=\{\infty\}\cup \{x_1,x_2,...,x_8\}\cup \{y_1,y_2,...,y_{14}\}$. Place a copy of an A2P maximum $8$-cycle packing of order $9$ on $\{\infty\}\cup \{x_1,x_2,...,x_8\}$ (by Example \ref {ex9}) and of order $15$ on $\{\infty\}\cup \{y_1,y_2,...,y_{14}\}$ (by Example \ref {ex15}) with the $4$-cycle leaves $(x_1,x_2,x_3,x_4)$ and $(y_1,y_2,y_3,y_4)$ and the $5$-cycle leaves $(y_9, y_{10}, y_{11}, y_{12}, y_{13})$. Then place an A2P maximum $8$-cycle packing of order $8$ on $\{x_1,x_2,x_3,x_4\} \cup \{y_1,y_2,y_3,y_4\}$ (by Example \ref {ex8}) with the $1$-factor leave $\{\{x_1,x_3\}, \{x_2,x_4\}, \{y_1,y_3\},\{y_2,y_4\}\}$, place an A2P $8$-cycle decomposition of $K_{4,4}$ on $\{x_1,x_2,$ $x_3,$  $x_4\} \cup \{y_5,y_6,y_7,y_8\}$, on $\{x_5,x_6,x_7,x_8\} \cup \{y_1,y_2,y_3,y_4\}$, and on $\{x_5,x_6,x_7,x_8\} \cup \{y_5,y_6,y_7,$ $y_8\}$, and finally place an A2P $8$-cycle decomposition of $K_{4,6}$ on $\{x_1,x_2,x_3,x_4\} \cup \{y_9,y_{10},...,y_{14}\}$ and on $\{x_5,x_6,x_7,x_8\} \cup \{y_9,y_{10},...,y_{14}\}$ by Lemma \ref{ts}.



\begin{lemma} \label{7mod16}
For every $n\equiv 7 \ ($mod $16)$ with $n \geq 23$, there exists an $A2P$ maximum $8$-cycle packing of order $n$.
\end{lemma}

\noindent {\bf Proof} There exist an A2P maximum $8$-cycle packing of order $23$ with a $5$-cycle leave by Example \ref{ex23}, and an A2P $8$-cycle decomposition of $K_{23}\setminus K_{7}$ by replacing $r=16$ and $s=6$ in Lemma \ref{hole}. Then the result follows by the Main Construction considering $h=7$. \hfill $\square$

\begin{center} \underline{\large $n\equiv 9$ (mod $16$)} \end{center}



\begin{lemma} \label{9mod16}
For every $n\equiv 9 \ ($mod $16)$, there exists an $A2P$ maximum $8$-cycle packing of order $n$.
\end{lemma}

\noindent {\bf Proof} There exist an A2P maximum $8$-cycle packing of order $9$ with a $4$-cycle leave by Example \ref{ex9} and an A2P $8$-cycle decomposition of $K_{25}\setminus K_{9}$  by Lemma \ref{hole} with $r=16$ and $s=8$. The result follows by Main Construction considering $h=9$. \hfill $\square$

\begin{center} \underline{\large $n\equiv 11$ (mod $16$)} \end{center}



\begin{lemma} \label{11mod16}
For every $n\equiv 11\ ($mod $16)$, there exists an $A2P$ maximum $8$-cycle packing of order $n$.
\end{lemma}

\noindent {\bf Proof} There exist an A2P maximum $8$-cycle packing of order $11$ with a $3$-cycle and a $4$-cycle leave by Example \ref{ex11} and an A2P $8$-cycle decomposition of $K_{27}\setminus K_{11}$ by  Lemma \ref{hole} with $r=16$ and $s=10$. The result follows by Main Construction considering $h=11$. \hfill $\square$

\begin{center} \underline{\large $n\equiv 13$ (mod $16$)} \end{center}

\begin{example} \label{ex13}
An $A2P$ maximum $8$-cycle packing of order $13$ exists.
\end{example}

\noindent To show this, let $\mathcal{X}=\{\infty\}\cup \{x_1,x_2,...,x_{8}\}\cup \{y_1,y_2,y_3,y_4\}$. Place an A2P maximum $8$-cycle packing of order $9$ (by Example \ref{ex9}) on $\{\infty\}\cup \{x_1,x_2,...,x_{8}\}$, with the $4$-cycle leave $(x_1,x_2,x_3,x_4)$. Then place an A2P maximum $8$-cycle packing of order $8$ (by by Example \ref{ex8}) on $\{x_1,x_2,x_3,x_4\} \cup \{y_1,y_2,y_3,y_4\}$ with the $1$-factor leave $\{\{x_1,x_3\}, \{x_2,x_4\}, \{y_1,y_3\}$, $\{y_2,y_4\}\}$. Then, place an A2P $8$-cycle decomposition of $K_{4,4}$ on $\{x_5,x_6,x_7,x_8\}\cup \{y_1,y_2,y_3,y_4\}$ by the Lemma \ref{K4t4s}. The unused edges $\{y_1,y_3\}$ and $\{y_2,y_4\}$ together with all the edges from $\infty$ to $\{y_1,y_2,y_3,y_4\}$ forms the leave bowtie $(y_1,y_3,\infty), (y_2,y_4,\infty)$.



\begin{lemma} \label{13mod16}
For every $n\equiv 13 \ ($mod $16)$, there exists an $A2P$ maximum $8$-cycle packing of order $n$.
\end{lemma}

\noindent {\bf Proof} There exist an A2P maximum $8$-cycle packing of order $13$ with a bowtie leave by Example \ref{ex13} and an A2P $8$-cycle decomposition of $K_{29}\setminus K_{13}$ by Lemma \ref{hole} used with $r=16$ and $s=12$. The result follows by Main Construction considering $h=13$. \hfill $\square$

\begin{center} \underline{\large $n\equiv 15$ (mod $16$)} \end{center}



\begin{lemma} \label{15mod16}
For every $n\equiv 15 \ ($mod $16)$, there exists an $A2P$ maximum $8$-cycle packing of order $n$.
\end{lemma}

\noindent {\bf Proof} There exist an A2P maximum $8$-cycle packing of order $15$ with a $4$-cycle and a $5$-cycle leave by Example \ref{ex15} and an A2P $8$-cycle decomposition of $K_{31}\setminus K_{15}$ by Lemma \ref{hole} $r=16$ and $s=14$. Then the result follows by the Main Construction considering $h=15$.\hfill $\square$

\begin{center} \underline{\large $n\equiv 0, 2,8$ and $10$ (mod $16$)} \end{center}

\begin{example} \label{ex10}
An $A2P$ maximum $8$-cycle packing of order $10$ exists.
\end{example}

\noindent Consider $\mathcal{X}=\mathbb{Z}_{10}$ and a maximum packing on $\mathcal{X}$ \\
$\mathcal{C}=\{(0,1,3,2,4,5,8,6)$, $(0,2,5,1,4,7,9,3)$, $(0,4,6,3,8,9,2,7)$, $(0,5,3,7,6,9,1,8)$, \\
$(1,6,2,8,4,9,5,7)\}$ with the leave $\mathcal{L}=\{\{1,2\},\{3,4\},\{5,6\},\{7,8\},\{9,0\}\}$ \\
Then one can choose the following inside cycles:\\
$\mathcal{C'}=\{(0,4,1,5,3,6,2,8)$, $(0,5,4,9,2,3,7,1)$, $(0,3,9,7,6,8,4,2)$,
$(0,7,5,9,8,3,1,6)$, \\ $(1,8,5,2,7,4,6,9)\}$.  \hfill $\square$

\begin{example} \label{ex16}
There exist $A2P$ maximum $8$-cycle packings of order $16$ and $18$.
\end{example}

\noindent There exist $A2P$ maximum $8$-cycle packings of order $8$ and $10$ with 1-factor leaves by Examples \ref{ex8} and \ref{ex10} respectively. There also exist  A2P $8$-cycle decomposition of $K_{8,8}$ and $K_{8,10}$ by Lemma \ref{K4t4s} and \ref{ts} respectively. Considering Lemma \ref{bipar} for $r=8$ with $s=8$ for the order $16$ and with $s=10$ for the order $18$ gives an A2P maximum $8$-cycle packing of orders $16$ and $18$ with a $1$-factor leave.



\begin{lemma} \label{02810}
For every $n\equiv 0,2,8$ and $10\ ($mod $16)$ with $n \geq 8$,  there exists an $A2P$ maximum $8$-cycle packing of order $n$.
\end{lemma}

\noindent {\bf Proof} There exist an A2P maximum $8$-cycle packing of orders $8,10, 16$ and $18$ by Examples \ref{ex8}, \ref{ex10}, and \ref{ex16}, respectively. An A2P maximum $8$-cycle packing of orders $24$ and $26$ exist by Lemma \ref{bipar} for $r=16$ with $s=8$ and $s=10$ respectively. For the same $r$ and $s$ in Lemma \ref{biparhole} we get $A2P$ maximum $8$-cycle decompositions of $K_{24}\setminus K_{8}$ and $K_{26}\setminus K_{10}$ respectively. Now we have all the ingredients to use the Main Construction for $h=0,2,8$ and $10$.\hfill $\square$

\begin{center} \underline{\large $n\equiv 4,6,12$ and $14$ (mod $16$)} \end{center}

\begin{example} \label{ex12}
There exist $A2P$ maximum $8$-cycle packings of orders $12, 14, 20$ and $22$.
\end{example}

\noindent Considering Lemma \ref{biparhole} for $r=8, 10, 14$ and $18$ with $s=4$, we get an A2P maximum $8$-cycle packings of $K_{12}\setminus K_{4}$,  $K_{14}\setminus K_{4}$, $K_{20}\setminus K_{4}$, and $K_{22}\setminus K_{4}$ with 1-factor leaves, respectively. But this is an A2P maximum $8$-cycle packing of orders $12, 14, 20$ and $22$, where the leave is a $K_4$ and a set of independent edges saturating the remaining elements.



\begin{lemma} \label{4,6,12,14}
For every $n\equiv 4,6,12$ and $14\ ($mod $16)$ with $n \geq 12$, there exists an $A2P$ maximum $8$-cycle packing of order $n$.
\end{lemma}

\noindent {\bf Proof} There exist an A2P maximum $8$-cycle packing of order $12,14, 20$, and $22$ by Example \ref{ex12} and an A2P maximum $8$-cycle packing of  $K_{20}\setminus K_{4}$, $K_{22}\setminus K_{6}$, $K_{28}\setminus K_{12}$ and $K_{30}\setminus K_{14}$ by Lemma \ref{biparhole} for  $r=16$ with $s=4, 6, 12,$ and $14$ respectively. Then, the result follows by the Main Construction by considering $h=4,6,12,14$.\hfill $\square$ \vspace{0.1in}

Now we can have the following result.

\begin{theorem}
There exists an $A2P$ maximum $8$-cycle packing of order $n$ for every $n\geq 8$.
\end{theorem}

\noindent {\bf Proof} Follows from the Lemma \ref{1mod16}, \ref{3mod16}, \ref{5mod16}, \ref{7mod16}, \ref{9mod16}, \ref{11mod16}, \ref{13mod16}, \ref{15mod16}, \ref{02810}, and \ref{4,6,12,14}.

\section{$8$-cycle packings that are not almost $2$-perfect}

Computer search shows that not all maximum packings with $8$-cycles are almost $2$-perfect. Even though all maximum packings with $8$-cycles for orders $8$ and $9$ are almost $2$-perfect, starting at order $10$, there are increasing number of maximum packings with $8$-cycles which do not carry this property. We generated about 2 million maximum packings with $8$-cycles of order 10, and only $0.35\%$ was not
almost $2$-perfect. On the other hand, computer search shows that when $n$ gets larger the percentage of $8$-cycle packings which are not almost $2$-perfect increases quite rapidly.

Below we give examples of $8$-cycle maximum packings which are not almost $2$-perfect for small orders. We will use them in the constructions to get $8$-cycle maximum packings which are not almost 2-perfect for all orders $n\geq10$.

\begin{example} \label{NA2P}
There exist $8$-cycle maximum packings which are not A$2$P of orders $10$, $11$, $12$, $13$, $15$, $16$ and $17$.
\end{example}

These packings of order $n$ are given on the set  $\{0,...,n-1\}$ for each $n$.

\noindent \textbf{Order 10}: $(0,2,1,3,4,6,5,7),(0,3,5,1,4,8,2,9)$, $(0,4,2,7,9,3,6,8)$, $(0,5,8,3,7,1,9,6)$,\\ $(1,6,2,5,9,4,7,8)$ with leave $\{\{0,1\},\{2,3\},\{4,5\},\{6,7\},\{8,9\}\}$.

\noindent\textbf{Order 11}: $(0,3,1,4,2,5,7,6)$, $(0,4,6,1,5,8,9,10)$, $(0,5,3,2,8,10,7,9)$, $(0,7,1,9,2,10,4,8)$, $(1,8,6,9,4,7,3,10)$, $(2,6,10,5,9,3,8,7)$ with leave $\{\{0,1\},\{1,2\}$, $\{0,2\},\{3,4\}$, $\{4,5\},\{5,6\}$, $\{3,6\}\}$.


\noindent\textbf{Order 12}: $(0,11,5,1,7,8,10,9),(0,4,6,2,9,11,3,10),(0,5,2,4,1,6,3,7),(0,6,5,3,4,7,2,8)$, $(1,8,3,9,4,10,2,11),(1,9,5,7,11,8,6,10),(4,8,5,10,7,9,6,11)$ with leave $K_4$ on $\{0,1,2,3\}$ and edges $\{\{4,5\},\{6,7\},\{8,9\},\{10,11\}\}$.

\noindent\textbf{Order 13}: $(0,11,5,1,7,8,10,9),(0,4,6,2,9,11,3,10),(0,2,4,1,3,5,6,7),(0,3,6,1,8,2,5,12)$, $(0,6,9,1,10,2,12,8),(1,11,2,7,3,8,4,12),(3,9,4,7,5,10,11,12)$, $(4,10,6,8,9,12,7,11)$, $(5,8$,\\$11$,$6,12,10,7,9)$ with leave $\{\{0,1\},\{1,2\},\{2,3\},\{3,4\},\{4,5\},\{0,5\}\}$.

\noindent\textbf{Order 15}: $(0,1,8,14,4,5,6,12)$, $(0,2,11,3,10,13,4,9))$, $(0,3,6,1,5,12,7,4)$, $(0,5,8,13,7,3,9,$ $6)$, $(0,8,2,6,13,9,1,10)$, $(0,11,14,12,2,9,5,13)$, $(9,12,3,1,4,10,2,14)$, $(1,13,2,4,8,10,5,14)$, $(0,7,1,11,4,6,10,14)$, $(1,2,3,5,7,11,8,12)$, $(2,5,11,13,3,8,6,7)$, $(3,4,12,10,11,9,7,14)$ with leave $\{\{7,8\},\{8,9\}$, $\{9,10\},\{7,10\},\{11,12\},\{12,13\},\{13,14\},\{14,6\},\{6,11\}\}$.

\noindent\textbf{Order 16}: Place a copy of the $8$-cycle packing of order $12$ given above on $\{0,...,11\}$ with the leave $K_4$ on $\{0,1,2,3\}$ and edges $\{\{4,5\},\{6,7\},\{8,9\},\{10,11\}\}$, and place a copy of any $8$-cycle packing of order $8$ on points $\{0,1,2,3,12,13,14,15\}$ with a 1-factor leave to cover the $K_4$ hole. Finally, place a copy of an $8$-cycle decomposition of $K_{8,4}$ on the bipartite graph with parts $\{4,5,6,7,8,9,10,11\}$ and $\{12,13,14,15\}$.

\noindent\textbf{Order 17}: The cyclic $8$-cycle system with base block $\{0,16,1,4,8,13,2,9\}$.

\begin{lemma} \label{lemNA2P}
If an $8$-cycle packing contains a subpacking which is not A$2$P, then the packing is also not A$2$P.
\end{lemma}
\noindent\textbf{Proof} Let $(\mathcal{X},\mathcal{C},\mathcal{L})$ be an A2P $8$-cycle packing with inside $8$-cycle packing $(\mathcal{X},\mathcal{C}',\mathcal{L})$. If $(\mathcal{X}_0,\mathcal{C}_0,\mathcal{L}_0)$ is a subpacking of $(\mathcal{X},\mathcal{C},\mathcal{L})$, then $(\mathcal{X}_0,\mathcal{C}'_0,\mathcal{L}_0)$ should be a subpacking of $(\mathcal{X},\mathcal{C}',\mathcal{L})$; as both $C_0$ and $C'_0$ are on the same vertex set $X_0$ and they both have the same number of $8$-cycles. So if $(\mathcal{X}_0,\mathcal{C}_0,\mathcal{L}_0)$ is not A2P then $(\mathcal{X},\mathcal{C},\mathcal{L})$ can not be A2P.\hfill $\square$

\begin{example} There exist $8$-cycle maximum packings of orders $19$ and $23 $ that are not A$2$P.
 \end{example}

 \noindent By replacing the packings of orders $11$ and $15$ in the constructions of Example \ref{ex19} and Example \ref{ex23} respectively with the packings of orders $11$ and $15$ given in Example \ref{NA2P}, one can construct packings that are not A2P for the orders $19$ and $23$.

 Observe that even if we let the edges in the leave to be used in the inside cycles (in other words let the leave changed in the formed packing), any collection of the inside $8$-cycles in packings of orders $11$ and $15$  in Example \ref{NA2P} can not form an $8$-cycle packing of orders $11$ or $15$.



\begin{theorem}
There exists an $8$-cycle packing which is not A$2$P for each $n\geq 10$.
\end{theorem}

\noindent\textbf{Proof} Considering the examples above for the constructions in Section 2 instead of the ones given in the previous section, one may construct an $8$-cycle maximum packing which is not A2P for each $n \geq 10$ by Lemma \ref{lemNA2P}.\hfill $\square$

\end{document}